\documentclass[11pt]{article}

\usepackage[utf8]{inputenc}
\usepackage[T1]{fontenc}
\usepackage{amsmath,amsfonts,amssymb,amsthm,bbm}
\usepackage{stackrel}
\usepackage{enumerate}
\usepackage{graphicx}
\usepackage{hyperref}
\usepackage{subfigure}
\usepackage{color}
\usepackage{bbm}

\newtheorem{theorem}{Theorem}

\newtheorem{proposition}[theorem]{Proposition}
\newtheorem{lemma}[theorem]{Lemma}

\newtheorem{remark}[theorem]{Remark}

\numberwithin{theorem}{section}
\numberwithin{equation}{section}
\numberwithin{figure}{section}

\newcommand{\ind}{\mathbbm{1}}

\newcommand{\calC}{\mathcal{C}}

\newcommand{\ZZ}{\mathbb{Z}}
\newcommand{\RR}{\mathbb{R}}

\newcommand{\PP}{\mathbb{P}}
\newcommand{\EE}{\mathbb{E}}

\newcommand{\Ball}{B} 
\newcommand{\Ann}{A} 
\newcommand{\din}{\partial^{\textrm{in}}} 
\newcommand{\dout}{\partial^{\textrm{out}}} 
\newcommand{\cluster}{\mathcal{C}}
\newcommand{\Ch}{\mathcal{C}_H} 
\newcommand{\arm}{\mathcal{A}} 
\newcommand{\lra}{\leftrightarrow}

\begin{document}

\title{On the four-arm exponent\\ for 2D percolation at criticality}

\author{Jacob van den Berg\footnote{CWI and VU University Amsterdam; E-mail: \texttt{J.van.den.Berg@cwi.nl}.}, Pierre Nolin\footnote{City University of Hong Kong; E-mail: \texttt{bpmnolin@cityu.edu.hk}. Partially supported by a GRF grant from the Research Grants Council of the Hong Kong SAR (project CityU11306719).}}

\date{}

\maketitle

\begin{center}
\emph{This paper is dedicated to the memory of Vladas Sidoravicius,\\ whose enthusiasm and dynamism have been very stimulating to us. \\[6mm]}
\end{center}

\begin{abstract}

For two-dimensional percolation at criticality, we discuss the inequality $\alpha_4 > 1$ for the polychromatic four-arm exponent (and stronger versions, the strongest so far being $\alpha_4 \geq 1 + \frac{\alpha_2}{2}$, where $\alpha_2$ denotes the two-arm exponent).
We first briefly discuss five proofs (some of them implicit and not self-contained) from the literature. Then we observe that, by combining two of them, one gets a completely self-contained (and yet quite short) proof.

\bigskip

\end{abstract}

\section{Introduction} \label{sec:intro}

In this paper we focus on site percolation on the square lattice $(\mathbb{Z}^2,\mathbb{E}^2)$. The vertices of this lattice are the points in $\RR^2$ with integer coordinates, and the edges in $\EE^2$ connect all pairs of vertices $v, v' \in \ZZ^2$ with $\|v-v'\|_2 = 1$ ($\|.\|_2$ denoting the usual Euclidean norm). However, note that the results would also hold on any two-dimensional lattice with enough symmetries, such as the honeycomb lattice, and also for bond percolation.

We are interested in upper bounds for the probability that two disjoint clusters connect neighbors of the origin to distance $n$, i.e. in lower bounds on the corresponding exponent. This exponent is called two-arm exponent in \cite{Cerf_2015} (a paper concerning dimensions $\geq 2$), but in two dimensions it is the same as what is usually called four-arm exponent: two open arms, one for each of the two open clusters, separated by two closed arms (ensuring that these two clusters are indeed not connected by an open path). We denote the corresponding exponent by $\alpha_4$. In the particular case of site percolation on the triangular lattice, this exponent is known to be equal to $\frac{5}{4}$ \cite{Smirnov_Werner_2001}, and this is widely believed to hold for all ``nice'' two-dimensional lattices (for site percolation, as well as for bond percolation).

For the square lattice it has been known for quite some time that $\alpha_4 > 1$. This strict inequality is related to the so-called noise sensitivity
of certain percolation phenomena (see Sections \ref{sec:BKS} and \ref{sec:Sch_St}). This inequality (and stronger versions)
has an interesting history, due to the diversity of the problems where four-arm probabilities (and their analogs in higher dimensions) 
played, play, or might play, a role (for instance, the uniqueness of the infinite cluster and the famous conjecture that $\theta(p_c) = 0$ for every dimension).

The first paper from which a proof of $\alpha_4 > 1$ can be (implicitly) obtained is (as several authors have mentioned, but without giving details) Kesten's celebrated scaling relations paper \cite{Kesten_1987}. We discuss in some detail in Section \ref{sec:Kesten_scaling_rel} how to do this. 
This method is quite technical and assumes much percolation background. Readers without such background are advised to skip that section.

In Section \ref{sec:earlier_proofs} we discuss parts of four other papers in the literature which, sometimes implicitly, provide a proof (some of them of the stronger result $\alpha_4 \geq 1 + \frac{\alpha_2}{2}$). Those proofs avoid the heavy near-critical machinery from \cite{Kesten_1987}. However, in most of these papers the four-arm inequality came up as a by-product or a necessary ingredient, and the authors have not always strived for optimizing simplicity or length of the proof.
Several of the proofs use a concentration result (which for this inequality is not needed) and/or a so-called arm-separation result: a result by Kesten which, although intuitively appealing, has a rather long and cumbersome proof.

A natural question is whether there is a short and self-contained proof that can be given in the first part of an introductory course on percolation theory, right after presenting the classical Russo-Seymour-Welsh result on crossing probabilities. We observed that one gets such a proof by following a special case of a proof by Garban in Appendix B of \cite{Schramm_Smirnov_2011} (which is inspired by a general inequality of \cite{ODonnell_Servedio_2007}, see also \cite{Garban_Steif_2015}), with modifications and ingredients from Cerf's arguments in \cite{Cerf_2015}.
This proof is presented in Section \ref{sec:proof}. It gives the stronger version of the inequality mentioned above, as stated more precisely in Theorem \ref{thm:main} below, but it is probably also, essentially, the shortest self-contained proof of the weaker version $\alpha_4 > 1$.

\begin{theorem} \label{thm:main}
For site percolation  on the square lattice $(\mathbb{Z}^2,\mathbb{E}^2)$ at criticality ($p = p_c^{\textrm{site}}(\mathbb{Z}^2)$), the following inequality between the two- and four-arm exponents, denoted by (resp.) $\alpha_2$ and $\alpha_4$, holds:
\begin{equation} \label{eq:main}
\alpha_4 \geq 1 + \frac{\alpha_2}{2}.
\end{equation}
\end{theorem}

We want to stress again that Theorem \ref{thm:main} is not new, but that the proof presented in Section \ref{sec:proof} (a modification and combination of other proofs) is arguably the most self-contained. It does not use Kesten's arm-separation results \cite{Kesten_1987}: in fact, it only uses pre-1980 percolation, namely the Russo-Seymour-Welsh result that at criticality, ``box-crossing probabilities are bounded away from $0$ and $1$''.

\subsection*{Organization of the paper}
In Section \ref{sec:percolation}, we set notation, and we recall the properties of critical percolation in 2D that we are going to use. We then comment on earlier (explicit or implicit) proofs of the inequality $\alpha_4 > 1$ (or even of \eqref{eq:main}) in Sections \ref{sec:Kesten_scaling_rel} and \ref{sec:earlier_proofs}, before turning to the 
self-contained proof of Theorem~\ref{thm:main} in Section~\ref{sec:proof}.

\section{Two-dimensional percolation at criticality} \label{sec:percolation}

\subsection{Setting and notations} \label{sec:setting}

\begin{figure}
\begin{center}

\subfigure{\includegraphics[width=.46\textwidth]{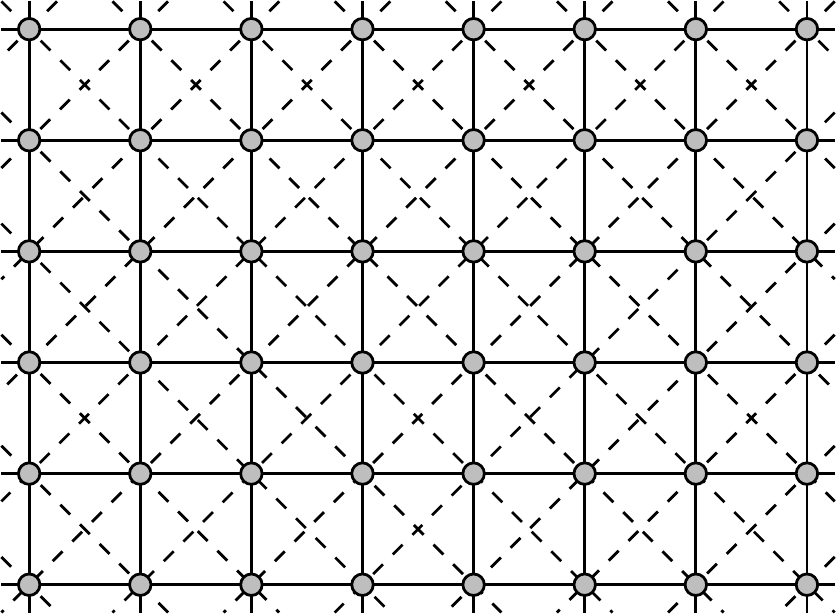}}
\hspace{0.04\textwidth}
\subfigure{\includegraphics[width=.46\textwidth]{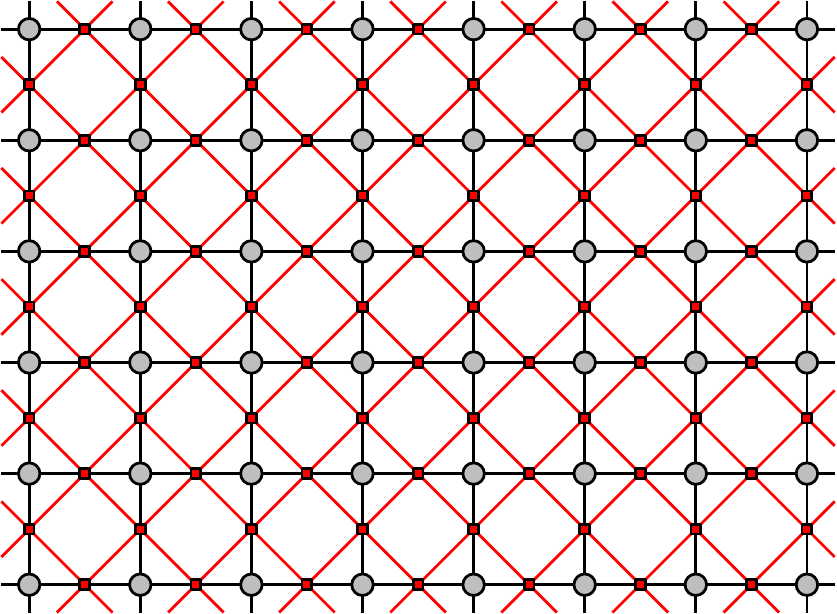}}\\

\caption{\label{fig:square_lattice} \emph{Left:} This figure shows the square lattice $G$, as well as the *-lattice obtained by adding the two diagonal edges (in dashed line) to every face of $G$. \emph{Right:} This figure depicts, in red, the medial lattice of $G$.}
\end{center}
\end{figure}

Recall that we work with the square lattice $G=(V,E)$, with set of vertices $V = \ZZ^2$, and set of edges $E = \EE^2$ connecting any two vertices which are at a Euclidean distance $1$ apart (i.e. differing along exactly one coordinate, by $\pm 1$). Two vertices $v, v' \in \ZZ^2$ are adjacent (or neighbors) if they are connected by an edge, i.e. $\{v,v'\} \in E$, and we write it $v \sim v'$. For a subset of vertices $A \subseteq V$, its inner and outer vertex boundaries are defined as, respectively,
$$\din A := \big\{ v \in A \: : \: v \sim v' \text{ for some } v' \in V \setminus A \big\}$$
and $\dout A := \din(V \setminus A)$. The matching lattice $G^*=(V^*,E^*)$, or simply *-lattice, is obtained from $G$ by adding the two diagonal edges to each face, as shown on Figure \ref{fig:square_lattice} (Left), and we use the notation $\sim^*$ for adjacency on $G^*$. A path (resp *-path) of length $k \geq 1$ on $G$ (resp. $G^*$) is a finite sequence of vertices $v_0, v_1, \ldots, v_k$ such that $v_i \sim v_{i+1}$ (resp. $v_i \sim^* v_{i+1}$) for all $i=0, \ldots, k-1$. We denote by $\Ball_n := [-n,n]^2$ the ball of radius $n \geq 0$ around $0$ for the $L^{\infty}$ norm $\|.\| = \|.\|_{\infty}$, and by $\Ann_{n_1, n_2} := \Ball_{n_2} \setminus \Ball_{n_1-1}$ the annulus with radii $0 \leq n_1 < n_2$ centered at $0$.

We also introduce the medial lattice $G^{\diamond} = (V^{\diamond},E^{\diamond})$ of $G$, for which a vertex $e^{\diamond} \in V^{\diamond}$ is located at the middle of every edge $e \in E$, and two such vertices $e^{\diamond}$, $e'^{\diamond}$ in $V^{\diamond}$ are connected by an edge \emph{if and only if} the corresponding edges $e$, $e'$ are incident to a common vertex in $V$: see Figure \ref{fig:square_lattice} (Right).

Bernoulli site percolation on $G$ with parameter $p \in [0,1]$ is obtained by declaring each vertex $v \in V$ either open or closed, with respective probabilities $p$ and $1-p$, independently of the other vertices. We denote by $\Omega := \{0,1\}^V$ the set of configurations $(\omega_v)_{v \in V}$, where $\omega_v = 1$ if $v$ is open, and $\omega_v = 0$ if $v$ is closed. We write $\PP_p$ for the product measure with parameter $p$ on $\Omega$.

Two vertices $v, v' \in V$ are connected (resp. *-connected) if there exists a path (resp. *-path) of length $k$, for some $k \geq 1$, along which all vertices are open (resp. closed), and we use the notation $v \lra v'$ (resp. $v \lra^* v'$). More generally for $A, A' \subseteq V$, $A \lra A'$ (resp. $A \lra^* A'$) means that there exist $v \in A$ and $v' \in A'$ such that $v \lra v'$ (resp. $v \lra^* v'$). Open vertices can be grouped into maximal connected components, that we call open clusters, and we denote by $\cluster(v)$ the open cluster containing a given $v \in V$ (with $\cluster(v) = \emptyset$ if $v$ is closed). Closed *-clusters are defined in a similar way.

\begin{figure}
\begin{center}

\includegraphics[width=.5\textwidth]{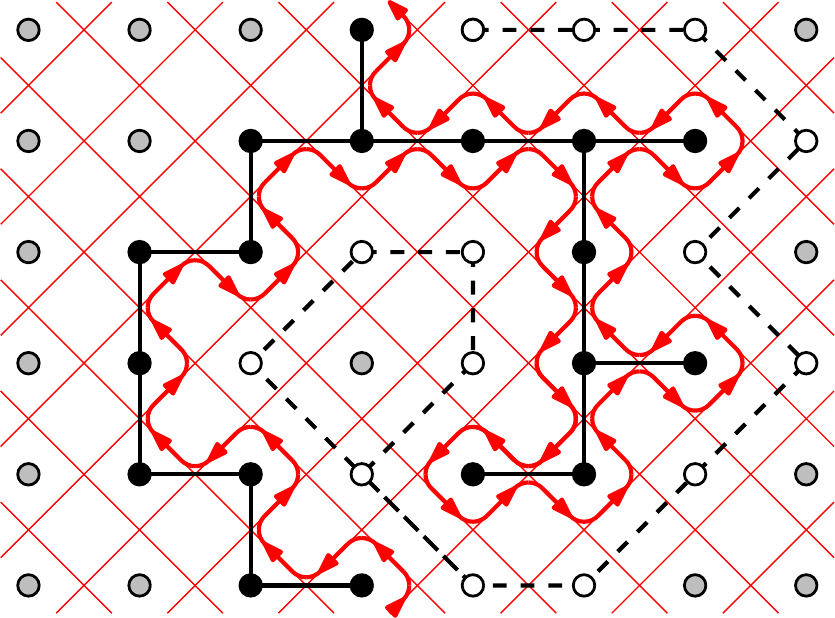}
\caption{\label{fig:exploration_process} An exploration process on $G$ following an ``interface'' between open and closed sites. It can be seen as an edge-self-avoiding path on the medial graph $G^{\diamond}$ of $G$ (each edge of $G^{\diamond}$ is followed at most once, although some vertices may be visited several times). The black and white vertices are revealed during the exploration, and respectively open and closed, while the grey vertices are left unexplored.}

\end{center}
\end{figure}

Exploration processes turn out to be an important ingredient in the proofs below. Such processes determine the outer boundary of an open cluster by revealing it in a step-by-step manner: all the open vertices along it, together with all the adjacent closed vertices (and discovering no other vertices). As shown on Figure \ref{fig:exploration_process}, they can be seen as edge-self-avoiding paths on the medial lattice $G^{\diamond}$.

Site percolation of $G$ displays a phase transition at a percolation threshold $p_c = p_c^{\textrm{site}}(G)$: for all $p < p_c$ there exists almost surely (a.s.) no infinite open cluster and a unique infinite closed *-cluster, while for all $p > p_c$ there is a.s. a unique infinite open cluster but no infinite closed *-cluster. In the present paper, we are concerned with the critical regime $p = p_c$, where neither infinite open clusters nor infinite closed *-clusters do exist. We refer the reader to the classical references \cite{Kesten_1982, Grimmett_1999} for more background on percolation theory.

Finally, the cardinality of a set $S$ is denoted by $|S|$, and for an event $E$, its indicator function $\ind_E$ is defined by: $\ind_E(\omega) = 1$ if $\omega \in E$, and $\ind_E(\omega) = 0$ otherwise.

\subsection{Critical regime} \label{sec:critical_regime}

We now recall classical definitions and properties concerning Bernoulli percolation at the critical point $p_c$.

If $R = [x_1,x_2] \times [y_1,y_2]$ (for some integers $x_1 < x_2$, $y_1 < y_2$) is a rectangle on the lattice, we denote by $\Ch(R)$ (resp. $\Ch^*(R)$) the existence of an open path (resp. closed *-path) in $R$ connecting the left side $\{x_1\} \times [y_1,y_2]$ and the right side $\{x_2\} \times [y_1,y_2]$. The classical Russo-Seymour-Welsh (RSW) theory states that
\begin{equation} \label{eq:RSW}
\PP_{p_c} \big( \Ch( [0,4n] \times [0,n] ) \big) \geq \delta_4 \quad \text{and} \quad \PP_{p_c} \big( \Ch^*( [0,4n] \times [0,n] ) \big) \geq \delta_4
\end{equation}
for some universal $\delta_4 > 0$. Using standard arguments, \eqref{eq:RSW} implies that for $\delta' = (\delta_4)^4 > 0$,
\begin{equation} \label{eq:RSW_csq}
\PP_{p_c}(\Ball_n \lra \din \Ball_{2n}) \leq 1 - \delta'.
\end{equation}
For $1 \leq n_1 < n_2$, let $\calC_{n_1,n_2}$ denote the collection of open clusters in $\Ball_{n_2}$ connecting $\Ball_{n_1}$ and $\din \Ball_{n_2}$. For future reference, observe that for some universal $c_1 < \infty$:
\begin{equation} \label{eq:upper_bd_C}
\text{for all } n \geq 1, \: \ell \geq n, \quad \EE_{p_c}\big[ |\calC_{n,n+\ell}|^2 \big] \leq c_1.
\end{equation}
Indeed, we know from \eqref{eq:RSW_csq} that $\PP_{p_c}(|\calC_{n,n+\ell}| \geq 1)$ is bounded away from $0$ and $1$, uniformly in $n$ and $\ell \geq n$. Hence, by the BK inequality, $|\calC_{n,n+\ell}|$ is (uniformly in $n$ and $\ell \geq n$) stochastically dominated by a geometrically distributed random variable, which gives \eqref{eq:upper_bd_C}.

Let $k \geq 1$, we consider the alternating sequence $\sigma_k = (oco\ldots) \in  \{o,c\}^k$, where $o$ and $c$ stand for ``open'' and ``closed'' respectively. In an annulus $A = \Ann_{n_1,n_2}$ ($0 \leq n_1 < n_2$), let $\arm_k(A)$ be the event that there exist $k$ disjoint paths $(\gamma_i)_{1 \leq i \leq k}$ in $A$, in counter-clockwise order, each connecting two vertices $v$ and $v'$ with $\|v\|  = n_1$ and $\|v'\|  = n_2$, and with respective types prescribed by $\sigma_k$ (i.e. $\gamma_i$ is an open path if $i$ is odd, and a closed *-path if $i$ is even). We write
\begin{equation} \label{eq:def_pi}
\pi_k(n_1,n_2) := \PP_{p_c}\big( \arm_k(\Ann_{n_1,n_2}) \big),
\end{equation}
and in particular $\pi_k(n) := \pi_k(\tilde{k},n)$, where $\tilde{k}$ is the smallest integer for which $\din \Ball_{\tilde{k}} \geq k$. Note that in this paper we consider only the cases $k=1,2,4$, for which $\tilde{k} = 0,1,1$ respectively. Finally, we introduce the $k$-arm (polychromatic, unless $k=1$) exponent
\begin{equation} \label{eq:def_exp}
\alpha_k := - \limsup_{n \to \infty} \frac{\log \pi_k(n)}{\log n}.
\end{equation}
It follows from standard constructions again (based on \eqref{eq:RSW}) that
$$\text{for all } k \geq 1, \quad \alpha_k \in (0,\infty).$$

\begin{remark} \label{rem:exponents}
\begin{enumerate}[(a)]
\item These arm exponents are known rigorously in the particular case of site percolation on the triangular lattice: $\alpha_1 = \frac{5}{48}$ \cite{Lawler_Schramm_Werner_2002}, and for all $k \geq 2$, $\alpha_k = \frac{k^2-1}{12}$ \cite{Smirnov_Werner_2001}. It is widely believed that these exponents should have the same values on other two-dimensional lattices such as the square lattice, considered in this paper.

\item \label{rem:separation} Adding certain ``macroscopic'' restrictions concerning the endpoints of the arms (for instance, in the case of four arms, that one endpoint is on the ``north'' side of $B_n$, and one on the west, one on the south, and one on the east side) does not increase the corresponding exponent. This ``arm-separation result'' was an important technical intermediate result by Kesten in his paper on scaling relations \cite{Kesten_1987}. Its proof is quite long and far from easy.
\end{enumerate}
\end{remark}

\section{Proof from Kesten's scaling relations (1987)} \label{sec:Kesten_scaling_rel}

In this section, we point out how the inequality $\alpha_4 > 1$ can be extracted from the results of \cite{Kesten_1987}. To the best of our knowledge, this paper is where the inequality $\alpha_4 > 1$ was first (implicitly) proved. Note that in this part, we assume much more percolation knowledge than in the rest of our paper, and the explanation below is mainly meant for specialists.

Other authors have already observed that the inequality $\alpha_4 > 1$ (or even better bounds on $\alpha_4$) can be obtained from \cite{Kesten_1987}. For instance, the paper \cite{Benjamini_Kalai_Schramm_1999} (that we discuss in more detail below, see Section \ref{sec:BKS}) says in Remark 4.2: ``Although this is better than the general bound \ldots, a somewhat better bound can be extracted from Kesten's\ldots''. But as far as we know, the authors did not write details about \emph{how} to obtain it from \cite{Kesten_1987}.

At first sight, doing so requires the assumption that some exponents exist. More explicitly, we assume first the existence of $\alpha_1$ (i.e. that the limit superior in \eqref{eq:def_exp} can be replaced by an actual limit), which implies that there is $\delta > 0$ such that
$$\PP_{p_c}(|\cluster(0)| \geq n) = n^{- \frac{1}{\delta} + o(1)} \quad \text{as $n \to \infty$}.$$
In addition, we need to assume the existence of $\alpha_4$, or equivalently of $\nu > 0$ such that $L(p) = |p-p_c|^{- \nu + o(1)}$ as $p \to p_c$, where the characteristic length $L$ is defined by $L(p) := \min \big\{ n \geq 1 \: : \: \PP_p \big( \Ch( [0,n] \times [0,n] ) \big) \leq 0.001 \big\}$ (resp. $\geq 0.999$) for $p < p_c$ (resp. $p > p_c$).

Corollary 2 in \cite{Kesten_1987} then states the inequality $\nu \geq \frac{\delta+1}{\delta}$. This inequality follows from previous results in \cite{Kesten_1987}, combined with either of the following two inequalities, as $p \nearrow p_c$:
\begin{equation} \label{eq:Durrett_Nguyen}
\frac{\mathbb{E}_p \big[ | \cluster(0) |^2 \big]}{\mathbb{E}_p \big[ | \cluster(0) | \big]} \geq (p_c - p)^{- 2 + o(1)}
\end{equation}
(see (3) in \cite{Durrett_Nguyen_1985}, Section 5), or
\begin{equation} \label{eq:Newman}
\mathbb{E}_p \big[ | \cluster(0) | \big] \geq (p_c-p)^{- 2 (\delta - 1) / \delta + o(1)}
\end{equation}
(see \cite{Newman_1987}, Theorem 1.3). Note that in \cite{Kesten_1987}, these inequalities \eqref{eq:Durrett_Nguyen} and \eqref{eq:Newman} are stated in terms of the critical exponents corresponding to the quantities in their l.h.s., usually denoted by $\Delta_2$ and $\gamma$ (respectively).

Hence, we have in particular $\nu > 1$. From the scaling relation $(2 - \alpha_4) \nu = 1$ (which follows from (4.28) and (4.33) in \cite{Kesten_1987}), we can thus obtain $2 - \alpha_4 < 1$, so the desired inequality $\alpha_4 > 1$. Moreover, we can actually get $\alpha_4 \geq 1 + \frac{\alpha_1}{2}$, by following more closely the previous sequence of inequalities and using the relation $\frac{2}{\delta+1} = \alpha_1$, proved in \cite{Kesten_1987b} (see the two sentences below (1.20) in \cite{Kesten_1987}, and note that in the notations of this paper, $\frac{1}{\delta_r}$ refers to the exponent $\alpha_1$).

Even if we do not assume the existence of some exponents, a large part of the results in \cite{Kesten_1987} can still be stated and established. In particular, one has the scaling relation
\begin{equation}
|p-p_c| L(p)^2 \pi_4(L(p)) \asymp 1
\end{equation}
as $p \to p_c$ (see (4.28) and (4.33) in \cite{Kesten_1987}, or Proposition 34 in \cite{Nolin_2008}). However, after closer inspection it is not immediately clear how to obtain the inequality $\alpha_4 \geq 1 + \frac{\alpha_1}{2}$ (or even $\alpha_4 > 1$).

We now explain how to obtain this inequality from the proof of \eqref{eq:Durrett_Nguyen} in \cite{Durrett_Nguyen_1985}. Note that if we try to follow the proof of \eqref{eq:Newman} in \cite{Newman_1987} instead, a difficulty arises. Indeed, the hypothesis (1.17) of Theorem 1.3 in \cite{Newman_1987} amounts to a lower bound on $\PP_p(|\cluster(0)| \geq n)$, while our definition of $\alpha_1$ involves an upper bound. As a consequence, we could not see how to use the reasonings in this paper (although it may be possible, we have not tried very hard).

Even though the paper \cite{Durrett_Nguyen_1985} (see Section 5) assumes the existence of exponents, we were able to fix this issue, and we now sketch briefly how to do it. For that, we use the (now-classical) scaling relations
\begin{equation} \label{eq:scaling_relations}
\chi(p) = \mathbb{E}_p \big[ | \cluster(0) | \big] \asymp L(p)^2 \pi_1(L(p))^2 \quad \text{and} \quad \mathbb{E}_p \big[ | \cluster(0) |^2 \big] \asymp L(p)^4 \pi_1(L(p))^3
\end{equation}
as $p \nearrow p_c$ (this is (1.25) in \cite{Kesten_1987}, for $t=1$ and $t=2$ respectively). In addition, one also has
\begin{equation}
\frac{d \chi(p)}{dp} \asymp L(p)^2 \pi_4(L(p)) \cdot \chi(p).
\end{equation}
Indeed, this can be proved by estimating $\frac{d}{dp} \PP_p(0 \lra v)$ for each $v \in \ZZ^2$, and then using similar reasonings as in \cite{Kesten_1987}. For $p < p_c$, these relations can be combined with the following inequality from \cite{Durrett_Nguyen_1985} (see p.266):
\begin{equation}
\mathbb{E}_p \big[ | \cluster(0) |^2 \big] \geq \frac{K}{\chi(p)} \bigg( \frac{d \chi(p)}{dp} \bigg)^2,
\end{equation}
for some universal constant $K \in (0,\infty)$. Hence, we get
\begin{equation}
\pi_4(L(p)) \leq K^{-1/2} L(p)^{-1} \pi_1(L(p))^{1/2}.
\end{equation}
Since $L(p) \to \infty$ as $p \nearrow p_c$, this gives the desired inequality between $\alpha_1$ and $\alpha_4$.

As a conclusion, we want to stress that one drawback of this approach is that it requires the arm-separation result mentioned in Remark \ref{rem:exponents}(\ref{rem:separation}). Also, we used quite heavy results on the behavior of percolation near criticality to deduce an inequality which is purely about the behavior at criticality. Proofs ``staying at criticality'' are arguably more satisfying.

\section{Other proofs in the literature} \label{sec:earlier_proofs}

We now discuss four papers in the literature which show lower bounds on $\alpha_4$ without using the quite heavy near-critical results in Kesten's paper \cite{Kesten_1987}.

The first three papers do this for bond percolation on the square lattice, and they are related to questions of noise sensitivity for a configuration at criticality. Presumably, after small modifications they also work for site percolation. We keep using the same notation $\pi_4(n)$ etcetera as we did for site percolation. These papers are: a paper by Benjamini, Kalai and Schramm \cite{Benjamini_Kalai_Schramm_1999} (Section~\ref{sec:BKS}), a paper by Schramm and Steif \cite{Schramm_Steif_2010} (Section~\ref{sec:Sch_St}), and an appendix by Garban in a paper by Schramm and Smirnov \cite{Schramm_Smirnov_2011} (Section~\ref{sec:Garban}). For some of the results in these papers, we also refer the reader to Sections~6.2.2 and 8.5 in the book \cite{Garban_Steif_2015} by Garban and Steif.

Finally, we discuss a paper by Cerf \cite{Cerf_2015} (Section~\ref{sec:Cerf}), which is written for site percolation on the square lattice (and, more generally, on the hypercubic lattice $\ZZ^d$ in any $d \geq 2$). Contrary to the above-mentioned papers, this paper is mostly concerned with dimensions $d \geq 3$, but, as we explain, it still yields interesting properties in dimension $d=2$.

Each of these papers uses some kind of exploration procedure in its proof of $\alpha_4 > 1$. And each of the first three papers uses Kesten's arm-separation result (see Remark \ref{rem:exponents}(\ref{rem:separation})). The proofs from \cite{Benjamini_Kalai_Schramm_1999} and \cite{Cerf_2015} use a concentration inequality, but the proofs in \cite{Schramm_Steif_2010} and \cite{Schramm_Smirnov_2011} do not. The main contribution by Garban in \cite{Schramm_Smirnov_2011} is a multi-scale version of Theorem \ref{thm:main} (see Lemma~\ref{lem:garban} below).

The proofs in \cite{Schramm_Steif_2010} and \cite{Schramm_Smirnov_2011} seem to be, partly or indirectly, influenced by \cite{Benjamini_Kalai_Schramm_1999}, but none of these three papers appears to be influenced by \cite{Aizenman_Kesten_Newman_1987} or \cite{Gandolfi_Grimmett_Russo_1988}. On the other hand, \cite{Cerf_2015} is influenced from these last two papers, but it seems to be completely independent of \cite{Benjamini_Kalai_Schramm_1999, Schramm_Steif_2010, Schramm_Smirnov_2011}.

Throughout this section the percolation parameter is equal to the bond or site (depending on the context) percolation threshold on the square lattice, and we omit it from our notation.

\subsection{The Benjamini-Kalai-Schramm paper (1999)} \label{sec:BKS}

The paper \cite{Benjamini_Kalai_Schramm_1999} is the first to give (for bond percolation on the square lattice) a proof of $\alpha_4 > 1$ without using the near-critical percolation results of \cite{Kesten_1987}.

Consider the event $A = A_m = \Ch( [0,m+1] \times [0,m] )$, and recall that an edge $e$ is said to be pivotal for $A$ if changing the state of $e$ changes the occurrence, or not, of $A$. The following is shown in \cite{Benjamini_Kalai_Schramm_1999}, where the only percolation knowledge used in the proof is the classical consequence from RSW that there exist $\rho, C > 0$ such that:
\begin{equation} \label{eq:RSW_csq2}
\text{for all } n \geq 1, \quad \PP_{p_c}(0 \lra \partial \Ball_n) \leq C n^{-1/\rho}
\end{equation}
(which follows immediately from \eqref{eq:RSW_csq}).

\begin{proposition}[\cite{Benjamini_Kalai_Schramm_1999}, equation (4.2) and Remark 4.2] \label{prop:crossing_BKS}
There is a constant $C > 0$ such that: for all $m \geq 1$,
\begin{equation} \label{eq:crossing_bks}
I(A) \leq C m^{1 - 1/{3 \rho}} (\log m)^{3/2},
\end{equation}
where $I(A)$ is the expected number of pivotal edges for the event $A$.
\end{proposition}

It follows from Kesten's arm-separation result that each edge in, say, the $\frac{m}{2} \times \frac{m}{2}$ square centered in the middle of the large box has a probability of order $\pi_4(m)$ to be pivotal. Since the expected number of pivotal edges in that square is smaller than or equal to the l.h.s. of \eqref{eq:crossing_bks}, we get $m^2 \pi_4(m) \leq C' m^{1 - 1/{3 \rho}} (\log m)^{3/2}$ (for some constant $C'$) and hence,
\begin{equation}
\pi_4(m) \leq C' m^{-1 - 1/{3 \rho}} (\log m)^{3/2}.
\end{equation}
Recalling the meaning of $\rho$, this gives, in our earlier notation,
\begin{equation}
\alpha_4 \geq 1 + \frac{\alpha_1}{3}.
\end{equation}

Proposition \ref{prop:crossing_BKS} is used in \cite{Benjamini_Kalai_Schramm_1999} to show that these box-crossing events are noise sensitive. An event $E \subseteq \Omega := \{0,1\}^n$ is said to be noise-sensitive if, roughly speaking, the following holds. For a large fraction of the configurations $\omega \in \Omega$, knowing $\omega$ does not significantly help to predict whether a perturbed configuration $\omega'$ (obtained from $\omega$ by randomly and independently flipping with small probability the ``bits'' $\omega_i$, $i = 1, \ldots , n$) belongs to the event $E$.

The proof of Proposition \ref{prop:crossing_BKS} is somewhat spread over different locations in the paper. As indicated above, the main concern of the paper is noise sensitivity. The paper contains some theorems of an ``algebraic'' flavour (involving discrete Fourier analysis), which give, for a quite general setting (i.e. not specifically for percolation) sufficient conditions for noise sensitivity. This type of results, combined with Proposition \ref{prop:crossing_BKS}, is essential to conclude noise sensitivity of the box-crossing events, but it is not needed for the proof of Proposition \ref{prop:crossing_BKS} itself. This makes it a bit hard to locate precisely those ingredients in the paper needed for the proof of Proposition~\ref{prop:crossing_BKS} itself.

Another type of results in the paper is of a more probabilistic nature and gives, again in a quite general setting, upper bounds for the total influence, which can then be used to check if the earlier mentioned conditions for noise sensitivity hold. One of the latter results, used for the proof of Proposition~\ref{prop:crossing_BKS}, is the following Lemma~\ref{lem:majority_BKS}. Let us first explain the notation in that lemma.

As before, $\Omega = \{0,1\}^n$, and the probability distribution considered is the product distribution with parameter $\frac{1}{2}$ (i.e. the uniform distribution on $\Omega$). For a function $f: \Omega \rightarrow [0,1]$, and a subset $K$ of $\{1, \ldots, n\}$, the notation $I_K(f)$ is used for $\sum_{k \in K} I_k(f)$, where 
$$I_k(f) = \frac{1}{2^n} \sum_{\omega \in \Omega} \big| f(\omega) - f(\omega^{(k)}) \big|,$$
with $\omega^{(k)}$ the configuration obtained from $\omega$ by flipping $\omega_k$ (note that if $f$ is the indicator function of an event, then $I_k(f)$ is the probability that $k$ is pivotal for that event).

Finally, $M_K$ is the majority function for $K$, which takes the value $1$ if the family $(\omega_i)_{i \in K}$ has more $1$'s than $0$'s, the value $-1$ if it has more $0$'s than $1$'s, and the value $0$ otherwise.

\begin{lemma}[\cite{Benjamini_Kalai_Schramm_1999}, Corollary 3.2 and Theorem 3.1] \label{lem:majority_BKS}
Let $K \subseteq \{1, \ldots, n\}$, and $f : \Omega \rightarrow [0,1]$ be monotone. Then, for some universal constant $C$, 
\begin{equation}
I_K(f) \leq C \sqrt{|K|} \,  \EE \big[ f M_K \big] \bigg(1 + \sqrt{- \log{\EE \big[ f M_K \big]}} \bigg).
\end{equation}
\end{lemma}

The proof of Lemma \ref{lem:majority_BKS} is self-contained and not very long (about one page), but certainly not obvious: it is a clever and surprising combination of nice elementary observations and standard concentration-like inequalities. 

The other important ingredient in \cite{Benjamini_Kalai_Schramm_1999} for the proof of Proposition~\ref{prop:crossing_BKS} is the following. This ingredient is very specific to the percolation setting mentioned before. Consider the $(m+1) \times m$ box in Proposition~\ref{prop:crossing_BKS} and the crossing event $A$ there.

\begin{lemma}[\cite{Benjamini_Kalai_Schramm_1999}, two lines before equation (4.2)] \label{lem:IM_bound}
For each subset $K$ of the set of edges in the right half of the $(m+1) \times m$ box, 
\begin{equation}
\EE \big[ \ind_A M_K \big] \leq C\, m^{-1/3 \rho} \log m,
\end{equation}
where $C$ is some universal constant.
\end{lemma}

Before we say a few words about the proof of Lemma~\ref{lem:IM_bound}, let us first see how Proposition~\ref{prop:crossing_BKS} follows. Combining Lemma~\ref{lem:IM_bound} and Lemma~\ref{lem:majority_BKS} gives immediately
$$I_K(A) \leq C \sqrt{|K|} \, m^{-1/3 \rho} (\log m)^{3/2}$$
for each subset $K$ of the set of edges in the right half of the $(m+1) \times m$ box. By symmetry, it then also holds for every $K$ in the left half of the box, and hence (with $C$ replaced by $C \sqrt 2$) for every $K$. Taking for $K$ the set of all edges of the box gives Proposition \ref{prop:crossing_BKS}.

As to the proof of Lemma \ref{lem:IM_bound}, it is practically self-contained; the only percolation knowledge that it uses is \eqref{eq:RSW_csq2}. The main ingredients of the proof of Lemma \ref{lem:IM_bound} are an exploration argument (for the existence of a horizontal crossing in the box), and some necessary quantitative work, again (as in the proof of Lemma \ref{lem:majority_BKS}) including some concentration-like inequalities. The main idea in the proof is that, to detect whether or not there is a horizontal crossing, typically a very small portion of $K$ is inspected. Indeed, in a simple exploration procedure, starting on the left side of the box, only edges of which at least one endpoint is connected to the left side of the box are inspected. Since each edge $e$ of $K$ is at a distance $\geq m/2$ from the left side of the box, the probability that it is inspected is at most of order $m^{-1/\rho}$. Using this it is shown that, typically, the ``surplus'' of $0$'s or $1$'s on the part of $K$ inspected by the algorithm, is much smaller than that on the rest of $K$, and therefore is unlikely to be decisive for the value of $M_K$. The mentioned concentration-like inequalities are used to make this precise.

\subsection{Four-arm results in the Schramm-Steif paper (2010)} \label{sec:Sch_St}

The paper \cite{Schramm_Steif_2010} studies the set of times at which an infinite cluster appears in a critical dynamical 2D percolation model. Noise sensitivity plays an important role in that study.

Some intermediate key results in this paper are stated in terms of discrete Fourier analysis (w.r.t. the Fourier-Walsh expansion). One such result is Theorem 1.8 in the paper. Let $\Omega = \{0,1\}^n$ and let $f : \Omega \rightarrow \RR$ be a function. Theorem 1.8 gives, for each $k \leq n$, an upper bound for the sum of the squares $\hat f(S)^2$ of the Fourier coefficients, over $S \subseteq \{1, \ldots, n\}$ with $|S| = k$. In the case where $k=1$ and $f$ is the indicator function of an increasing event $A$, one can use (as mentioned in the remark below Theorem 4.1 in \cite{Schramm_Steif_2010}) that $\hat f(\{i\})$ is equal to the probability that $i$ is pivotal for $A$. For that special case, Theorem 1.8 in \cite{Schramm_Steif_2010} is as follows.

\begin{lemma}[special case of \cite{Schramm_Steif_2010}, Theorem 1.8] \label{lem:ScSt1.8}
Let $\Omega = \{0,1\}^n$ and let $E \subseteq \Omega$ be an increasing event. Further, let $A$ be a randomized algorithm which determines, by a step-by-step procedure, whether a configuration $\omega$ belongs to $E$ or not, and where at each step of the procedure, the value of exactly one $\omega_i$ is ``revealed'' (the choice of $i$ may depend on the values of the $\omega_j$'s that have already been inspected at that stage). The algorithm stops as soon as it is known whether $E$ occurs or not. Let $\delta_A$ be the maximum over all $i \in \{1, \ldots, n\}$ of the probability that $i$ is inspected. Then
\begin{equation} \label{eq:ScSt1.4}
\sum_{i=1}^n \PP \big( i \text{ is pivotal for } E \big)^2 \leq \delta_A \, \PP(E).
\end{equation}
\end{lemma}

The proof of Theorem 1.8 in \cite{Schramm_Steif_2010} is not long, and it is reasonably self-contained but quite subtle.

Another result in \cite{Schramm_Steif_2010} which is relevant for obtaining bounds on four-arm probabilities is Theorem 4.1 in that paper. It gives a suitable ``decision algorithm'' $A$ for the event that there is a horizontal open crossing of an $R \times R$ square. This algorithm needed special care because $\delta_A$ is the maximum revealment probability over \emph{all} edges in the square (not only the edges in the concentric $\frac{R}{3} \times \frac{R}{3}$ square). More precisely, Theorem 4.1 says (in our notation) the following.

\begin{lemma}[\cite{Schramm_Steif_2010}, Theorem 4.1] \label{thm:ScSt4.1}
For the above mentioned crossing event for site percolation on the triangular lattice, there is an algorithm $A$ with $\delta_A \leq R^{-\frac{1}{4} + o(1)}$. For the similar event for bond percolation on the square lattice, there exists a constant $a > 0$ and an algorithm $A$ with $\delta_A \leq R^{-a + o(1)}$.
\end{lemma} 

The paper \cite{Schramm_Steif_2010} gives a proof for the statement on the triangular lattice, and says that the proof of the statement for the square lattice is similar. Note that the value $\frac{1}{4}$ in Lemma \ref{thm:ScSt4.1} is the two-arm exponent $\alpha_2$ on the triangular lattice. From the proof of the lemma, it is not clear whether, in the case of the square lattice, we may take $a = \alpha_2$ in the above theorem. However, this \emph{is} clear for the weaker lemma where $\delta_A$ is replaced by the maximum revealment probability over the edges in the earlier mentioned $\frac{R}{3} \times \frac{R}{3}$ square. Combining that weaker lemma with a suitable modification of Lemma~\ref{lem:ScSt1.8} (where for $E$ we take the event that there is an open crossing of an $R \times R$ square, we replace the sum in the l.h.s. of \eqref{eq:ScSt1.4} by the smaller sum restricted to the vertices in the concentric $\frac{R}{3} \times \frac{R}{3}$ box, and $\delta_A$ is replaced as mentioned above), and then using Kesten's arm-separation result, gives $R^2 \, \pi_4(R)^2 \leq  R^{-\alpha_2 + o(1)}$, and hence Theorem \ref{thm:main}. See Corollary A.4 of \cite{Vanneuville_2019} for such modifications.

\subsection{The result of Garban (2011)} \label{sec:Garban}

In Appendix B of the paper \cite{Schramm_Smirnov_2011} by Schramm and Smirnov, Garban gives a ``multi-scale bound'' on the four-arm probability for bond percolation on $\ZZ^2$. More precisely, let $\varepsilon$ be such that there is a constant $c' > 0$ for which: for all $1 \leq m \leq n$, $\pi_2(m,n) \leq c' \big( \frac{m}{n} \big)^{2 \varepsilon}$. The following is proved in \cite{Schramm_Smirnov_2011}.
\begin{lemma}[\cite{Schramm_Smirnov_2011}, Appendix B] \label{lem:garban}
There is a constant $c > 0$ such that:
\begin{equation} \label{eq:multi_scale}
\text{for all } 1 \leq m \leq n, \quad \pi_4(m,n) \leq c \bigg( \frac{m}{n} \bigg)^{1 + \varepsilon}.
\end{equation} 
\end{lemma}

For the special case $m=1$, this gives $\alpha_4 \geq 1 + \frac{\alpha_2}{2}$. 
A nice aspect of Garban's proof is that it is completely focused on the problem in question, while the mentioned four-arm results in \cite{Benjamini_Kalai_Schramm_1999} and \cite{Schramm_Steif_2010} were in some sense (versions of) intermediate results needed in the proof of some other results. 

Interestingly, Garban says that the case $m=1$ ``can be extracted from \cite{Kesten_1987} as well as \cite{Benjamini_Kalai_Schramm_1999} or \cite{Schramm_Steif_2010}''. In fact, following his proof, but (roughly speaking) taking everywhere $m =1$, is considerably simpler than extracting a full proof for that case from the mentioned papers. Apart from the fact that it uses Kesten's arm-separation results, it is probably the shortest and most elegant proof that $\alpha_4 \geq 1 + \frac{\alpha_2}{2}$. It avoids concentration results (which were used in Cerf's computation, see the next section). As Garban indicates, a key part in his proof, in that special case $m=1$, is essentially an application of (or almost ``equivalent'' to the proof of) a quite general inequality of \cite{ODonnell_Servedio_2007} (see also the remark following the proof of Proposition~6.6 in Section~8.5 of \cite{Garban_Steif_2015}).

\subsection{A result by Cerf (2015)} \label{sec:Cerf}

Lemma 5.2 in \cite{Cerf_2015}, that we now state in any dimension $d \geq 2$, gives the following result (recall that $\calC_{n,n+\ell}$ is the collection of open clusters in $\Ball_{n + \ell}$ connecting $\Ball_n$ and $\din \Ball_{n + \ell}$).
\begin{lemma}[\cite{Cerf_2015}, Lemma 5.2] \label{lem:Cerf_5.2}
Let $d \geq 2$, and consider site percolation on the hypercubic lattice $\ZZ^d$. For all $p \in (0,1)$, $n \geq 1$ and $\ell \geq 0$,
\begin{align}
\PP_p\big( \arm_4( & \Ann_{1,2n+\ell}) \big) \nonumber \\[1mm]
& \leq \frac{2d (\log n)}{\sqrt{|\Ball_n|}} \EE_p \bigg[ \sqrt{|\calC_{n,n+\ell}|} \bigg] + \frac{4d}{p(1-p)} |\Ball_n|^2 e^{-2 (\log n)^2 p^2 (1-p)^2}. \label{eq:Cerf_5.2}
\end{align}
\end{lemma}
Note that this result holds for any $p \in (0,1)$. For our purpose, we will restrict, but only later, to $d=2$ and $p = p_c^{\textrm{site}}(\mathbb{Z}^2)$.

The proof of this lemma in \cite{Cerf_2015} is completely self-contained, it assumes no percolation knowledge at all.
It is a nice mixture of arguments with a combinatorial flavor, and application of a concentration inequality (see our comments later in this section).
As Cerf remarks, a version of this result, with only the parameter $n$, not $\ell$ (or, more precisely, with $\ell = 0$), is somewhat hidden
in the arguments of Gandolfi, Grimmett and Russo \cite{Gandolfi_Grimmett_Russo_1988} and Aizenman, Kesten and Newman \cite{Aizenman_Kesten_Newman_1987}, to 
prove the uniqueness of the infinite open cluster.

Following \cite{Cerf_2015}, taking $\ell = 0$ in \eqref{eq:Cerf_5.2} and using the trivial upper bound $| \din \Ball_n | \asymp n^{d-1}$ for $|\calC_{n,n+\ell}|$ gives
\begin{equation} \label{eq:Cerf_upper_bd1}
\PP_p\big( \arm_4( \Ann_{1,2n}) \big) \leq c \, \frac{\log n}{\sqrt{n}},
\end{equation}
where $c$ depends on the dimension $d$ only.

The main contribution in \cite{Cerf_2015} is to ``bootstrap'' \eqref{eq:Cerf_5.2} in a clever way: the inequality \eqref{eq:Cerf_upper_bd1} is used to improve the above-mentioned trivial upper bound for $\EE_p \big[ \sqrt{|\calC_{n,n+\ell}|} \big]$, which is then plugged into \eqref{eq:Cerf_5.2} to get an improvement of \eqref{eq:Cerf_upper_bd1},
then leading to an even better bound for $\EE_p \big[ \sqrt{|\calC_{n,n+\ell}|} \big]$, and so on. 
The introduction by Cerf of the extra parameter $\ell$ seems to provide the
flexibility needed to do this bootstrapping.

As pointed out in \cite{Cerf_2015}, for $d=2$ the final result obtained in this way is $\alpha_4 \geq \frac{11}{21}$, which looks disappointing.
However, the main focus in the paper is on dimensions $d \geq 3$, where the ``bootstrapping'' that we just explained does give interesting new results.

Nevertheless, it may be worth mentioning that, as we observed, \eqref{eq:Cerf_5.2} (and a modified version obtained from small changes in its proof) is also useful for the case $d=2$ (even without using the bootstrapping), as we point out now.

First, note that for $d=2$ and $p = p_c^{\textrm{site}}(\mathbb{Z}^2)$, $\EE_p \big[ \sqrt{|\calC_{n,2n}|} \big]$ is uniformly bounded in $n$ (so bootstrapping makes no sense for $d=2$). So, for $d=2$, \eqref{eq:Cerf_5.2}, now with $\ell = n$, actually gives
$$\pi_4(3n) \leq \tilde{c} \, \frac{\log n}{n}$$
for some constant $\tilde{c}$, and hence $\alpha_4 \geq 1$.

As we point out next, one can, with a very small modification in the proof of \eqref{eq:Cerf_5.2}, obtain $\alpha_4 \geq 1 + \frac{\alpha_1}{2}$. Lines 8--9 in Section 5 of \cite{Cerf_2015} give an upper bound for the quantity
\begin{equation} \label{eq:Cerf_sum1}
\sum_{C \in \calC} \sqrt{|\bar{C} \cap \Ball_n|},
\end{equation}
where $\calC = \calC_{n, 2n}$ and we denote $\bar{C} := C \cup \dout C$. Namely (by Jensen's inequality), this quantity is at most
\begin{equation} \label{eq:Cerf_sum2}
\sqrt{|\calC|} \sqrt{\sum_{C \in \calC} |\bar{C} \cap \Ball_n|},
\end{equation}
which, since every vertex $v$ belongs to at most $2d$ subsets $\bar{C}$ with $C \in \calC$, is at most $\sqrt{|\calC|} \sqrt{2d} \sqrt{|\Ball_n|}$. So for the expectation of the sum in \eqref{eq:Cerf_sum1}:
\begin{equation} \label{eq:Cerf_sum3}
\EE_{p_c} \Bigg[ \sum_{C \in \calC} \sqrt{|\bar{C} \cap \Ball_n|} \Bigg] \leq \EE_{p_c} \Big[ \sqrt{|\calC|} \Big] \sqrt{2d} \sqrt{|\Ball_n|},
\end{equation}
which is used later in \cite{Cerf_2015} to obtain \eqref{eq:Cerf_5.2}.

The ``very small modification'' that we meant is the following: by the Cauchy-Schwarz inequality, the expectation of \eqref{eq:Cerf_sum2} is at most
\begin{equation} \label{eq:Cerf_sum4}
\EE_{p_c} \Bigg[ \sqrt{|\calC|} \sqrt{\sum_{C \in \calC} \big| \bar{C} \cap \Ball_n \big|} \Bigg] \leq \sqrt{\EE_{p_c} \big[ |\calC| \big]} \sqrt{ \EE_{p_c} \Bigg[ \sum_{C \in \calC} \big| \bar{C} \cap \Ball_n \big| \Bigg]}.
\end{equation}
Since every $v \in \bigcup_{C \in \calC} (\bar{C} \cap \Ball_n)$ has an open path to $\din \Ball_{2n}$, the expectation in the second factor in \eqref{eq:Cerf_sum4} above is at most $2d |\Ball_n| \pi_1(n)$. So we get that the expectation of \eqref{eq:Cerf_sum1} is at most

\begin{equation} \label{eq:Cerf_sum5}
\EE_{p_c} \Bigg[ \sum_{C \in \calC} \sqrt{|\bar{C} \cap \Ball_n|} \Bigg] \leq \sqrt{\EE_{p_c} \big[ |\calC| \big]} \sqrt{2d} \sqrt{|\Ball_n|} \sqrt{\pi_1(n)}.
\end{equation}
Comparing this with the r.h.s. of \eqref{eq:Cerf_sum3} (and recalling that, for $d=2$, $\EE_{p_c} [ |\calC| ]$ is uniformly bounded), we see that we made appear an extra factor $\sqrt{\pi_1(n)}$. This then also causes the same additional factor in the first term in the r.h.s. of \eqref{eq:Cerf_5.2}, and yields
$$\alpha_4 \geq 1 + \frac{\alpha_1}{2}.$$

Finally, one gets (still for $d=2$) a further improvement by considering, in the proof in \cite{Cerf_2015}, instead of $\bar{C}$, the set of vertices $\tilde{C} := C^* \cup C'$, where
$$C^* := \big\{ v^* \in \dout C \cap \Ball_{2n} \: : \: v^* \lra^* \din \Ball_{2n} \big\},$$
and
$$C' := \big\{ v \in C \: : \: v \sim v^* \text{ for some } v^* \in C^* \big\},$$
and then using that from every $v \in \bigcup_{C \in \calC} (\tilde{C} \cap \Ball_n)$, one can find an open path and a closed *-path (starting from neighbors of $v$) to $\din \Ball_{2n}$. This now produces, instead of the above-mentioned $\sqrt{\pi_1(n)}$, an extra factor $\sqrt{\pi_2(n)}$ in the first term in the r.h.s. of \eqref{eq:Cerf_5.2},
so that we get $$\alpha_4 \geq 1 + \frac{\alpha_2}{2}.$$
Comparing the case $m=1$ of Garban's proof (mentioned in Section \ref{sec:Garban}) of this inequality with the proof in \cite{Cerf_2015} of \eqref{eq:Cerf_5.2},
we observe that the latter avoids Kesten's arm-separation result, and is thus more self-contained. It uses a large-deviation argument which makes it longer,
and which is, presumably, only useful for the case $d \geq 3$.

In the next section, we give a short and self-contained proof of $\alpha_4 \geq 1 + \frac{\alpha_2}{2}$, which can be considered as a combination of the proof of \eqref{eq:multi_scale} (in the special case $m=1$) in \cite{Schramm_Smirnov_2011} and the proof of \eqref{eq:Cerf_5.2} in \cite{Cerf_2015}.

\section{A self-contained proof of Theorem~\ref{thm:main}, based on Garban's and Cerf's arguments} \label{sec:proof}

\subsection{Introductory remarks}

We follow Garban's proof for the result in Section \ref{sec:Garban}, but restrict to the case $m=1$, 
and replace the \emph{event that there is a horizontal crossing of a box}, by the {\emph{number of connected components crossing an annulus}. 
The proof of Theorem \ref{thm:main} obtained in this way is, in some sense, a mixture of Garban's argument and that by Cerf:
it still exploits, as in Garban's proof (which, as said, was inspired by \cite{ODonnell_Servedio_2007}), the full power of symmetry provided by involving the notion of pivotality, while it also uses the 
advantage of considering the number of crossings of an annulus (as Cerf did) instead of the event (considered by Garban) that there is a horizontal crossing of a box. 
This enables one to avoid Kesten's arm-separation result (we do not see how to avoid that result in the proof of Lemma \ref{lem:garban} for a general $m \geq 1$).
To underline the flexibility of the method, we deal with \emph{site} percolation on the square lattice (which has less symmetry than bond percolation on that lattice), with parameter $p_c = p_c^{\textrm{site}}(\mathbb{Z}^2)$.

\subsection{Proof}
Let $n$ be a positive integer, and let $\Omega = \{0,1\}^{\Ball_{2 n}}$ be the set of all configurations of open and closed vertices in the box $\Ball_{2 n}$. Let $Z = |\calC_{n, 2n}|$ be the number of open clusters in $\Ball_{2 n}$ that have at least one vertex in each of $\Ball_n$ and $\din \Ball_{2 n}$. From \eqref{eq:upper_bd_C}, we know that for some universal $\bar c > 0$ (independent of $n$),
\begin{equation} \label{eq:G-EZ-bnd}
\EE_{p_c} \big[ Z^2 \big] \leq {\bar c}^2.
\end{equation}

Note that if we close an open vertex in $\Ball_{n-1}$, the value of $Z$ does not decrease. Let $v_1, v_2, \ldots$ be a list of the vertices in $\Ball_{n-1}$. For each $1 \leq j \leq |\Ball_{n-1}|$, define the random variable $C_j$ as follows:
$$
C_j = \begin{cases}
-(1 - p_c) & \text{if } v_j \text{ is open,} \\[1mm]
p_c & \text{if } v_j \text{ is closed.}
\end{cases}
$$

In the remainder of this proof, $\{v_j \text{ is pivotal}\}$ denotes the event that if the state of $v_j$ is changed, then the value of $Z$ changes as well. More precisely,
$$\big\{ v_j \text{ is pivotal} \big\} := \big\{ \omega \in \Omega \: : \: Z(\omega^{(j)}) \neq Z(\omega) \big\},$$
where $\omega^{(j)}$ denotes the configuration obtained from $\omega$ by ``flipping'' $\omega_{v_j}$.

We now consider an exploration procedure $\Gamma$ which counts the number $Z$ of open clusters in $\calC_{n, 2n}$. Roughly speaking, $\Gamma$ is constructed so as to follow successively the boundaries (as depicted on Figure \ref{fig:exploration_process}) of all open clusters in $\Ball_{2 n}$ that intersect $\din \Ball_{2n}$, starting from $\din \Ball_{2n}$. It has the property that each time it reaches a ``fresh'' vertex, the state of this vertex is revealed, open with probability $p_c$ and closed with probability $1-p_c$, independently of all information obtained so far in the procedure. We refer to Figure \ref{fig:exploration_gamma}, which shows an intermediate stage of this procedure, and where the vertices pivotal for $Z$ are marked.

\begin{figure}[h]
\begin{center}

\includegraphics[width=.7\textwidth]{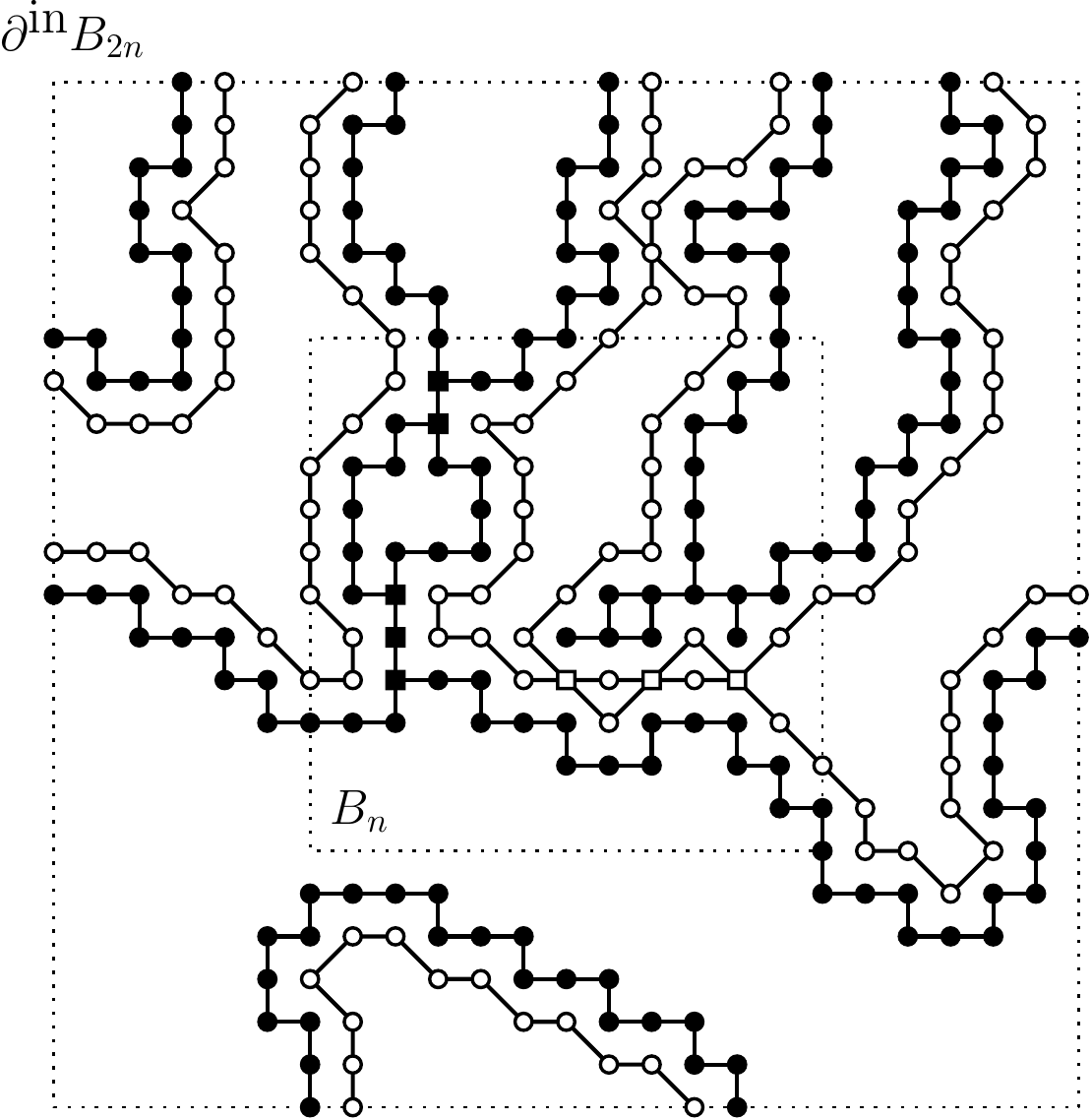}
\caption{\label{fig:exploration_gamma} This figures shows part of the exploration procedure $\Gamma$, which explores iteratively the ``interfaces'' between open clusters and closed *-clusters connected to $\din \Ball_{2n}$. The black vertices are open, and the white ones are closed. The vertices indicated with a square are pivotal for $Z$: changing the state of such a vertex $v$ would increase or decrease the value of $Z$, depending on whether $v$ is open or closed, respectively.}

\end{center}
\end{figure}

We let
$$Y_j := \ind_{v_j \text{ is visited by } \Gamma}.$$
Note that for each vertex $v$ visited by $\Gamma$ (and away from $\din \Ball_{2n}$), it is possible to find an open path and a closed *-path from neighbors (or *-neighbors) of $v$ to $\din \Ball_{2n}$. Since each vertex in $\Ball_n$ is at a distance at least $n$ from $\din \Ball_{2n}$, we obtain
\begin{equation} \label{eq:G-B6}
\EE_{p_c} \big[ Y_j \big] \leq c \, \pi_2(n)
\end{equation}
for some constant $c > 0$.

By the nature of the exploration path (the next step of the path depends only on the states of the vertices hit by the path so far),
\begin{equation}\label{eq:G-Y-equality}
Y_j(\omega^{(j)}) = Y_j(\omega).
\end{equation}
In particular, $C_j$ and $Y_j$ are independent, and $\EE_{p_c}[C_j Y_j] = 0$. For essentially the same reason, if $v_i$ and $v_j$ are two distinct vertices, then, at the first step in the procedure that one of these two vertices is hit, the $Y$- and $C$-values of the other vertex are conditionally independent, given all information obtained during the exploration so far. 
Because of this (and a similar argument for the case where neither $v_i$ nor $v_j$ is hit), we get:
\begin{equation} \label{eq:G-Bij}
\text{for all } i \neq j, \quad \EE_{p_c} \big[ (C_i Y_i) (C_j Y_j) \big] = 0.
\end{equation}

We now study $\EE_{p_c}[Z C_j Y_j]$ (this is analogous to Garban's proof, but with $Z$ instead of the indicator function of a crossing event). Clearly,
\begin{equation} \label{eq:G-B4a}
\EE_{p_c} \big[ Z C_j Y_j \big] = \EE_{p_c} \big[ Z C_j Y_j \, \ind_{v_j \text{ is pivotal}} \big] + \EE_{p_c} \big[ Z C_j Y_j \, \ind_{v_j \text{ is not pivotal}} \big].
\end{equation}
Let $\omega \in \Omega$. As is easy to check (using \eqref{eq:G-Y-equality}), we have 
$$\PP_{p_c}(\omega) C_j(\omega) Y_j(\omega) = - \PP_{p_c}(\omega^{(j)}) C_j(\omega^{(j)}) Y_j(\omega^{(j)}).$$

On the one hand, if $\omega \in \{v_j$ is not pivotal$\}$, then $\omega^{(j)} \in \{v_j$ is not pivotal$\}$ as well, and $Z(\omega) = Z(\omega^{(j)})$. Hence, the contribution of the pair $(\omega, \omega^{(j)})$ to the second term in the r.h.s. of \eqref{eq:G-B4a} is $0$, from which it follows that this term is equal to $0$. On the other hand, if $\omega \in \{v_j \text{ is pivotal}\}$, the state of $v_j$ must be explored by $\Gamma$. Hence, the first term of \eqref{eq:G-B4a} is equal to $\EE_{p_c} \big[ Z C_j \, \ind_{v_j \text{ is pivotal}} \big]$.

Now let $\omega \in \{v_j \text{ is pivotal}\}$, and suppose that $\omega_{v_j} = 1$, so that $C_j(\omega) = - (1-p_c)$. Then also $\omega^{(j)} \in \{v_j \text{ is pivotal}\}$, but $C_j(\omega^{(j)}) = p_c$. It follows that the contribution of the pair $(\omega, \omega^{(j)})$ to the first term in \eqref{eq:G-B4a} is $p_c (1-p_c) q \, (Z(\omega^{(j)}) - Z(\omega))$, where $q = q(\omega)$ denotes the probability of the configuration $(\omega_v)_{v \in \Ball_{2 n} \setminus \{v_j\}}$ (note that $q(\omega) = q(\omega^{(j)})$). Using that $Z(\omega^{(j)}) - Z(\omega) \geq 1$, and summing over all configurations in the event $\{v_j \text{ is pivotal}\}$, we obtain that the first term in the r.h.s. of \eqref{eq:G-B4a} is larger than or equal to
$$p_c (1-p_c) \PP_{p_c}(v_j \text{ is pivotal}).$$

By the above, and also observing that $\PP_{p_c}(v_j \text{ is pivotal}) \geq \pi_4(3 n)$ (indeed, if a vertex $v \in \Ball_n$ has four arms to distance $3n$, then it has four arms to $\din \Ball_{2n}$, and so it is pivotal for $Z$), we conclude that
\begin{equation} \label{eq:G-B8}
\EE_{p_c} \big[ Z C_j Y_j \big] = \EE_{p_c} \big[ Z C_j  \, \ind_{v_j \text{ is pivotal}} \big] \geq p_c (1 - p_c) \pi_4(3 n).
\end{equation}
The sum over $j$ of the l.h.s. of \eqref{eq:G-B8} satisfies (for some constant $\hat c > 0$):
\begin{align*}
\sum_j \EE_{p_c} \big[ Z C_j Y_j \big] & \leq \sqrt{ \EE_{p_c} \big[ Z^2 \big] \EE_{p_c} \Bigg[ \bigg( \sum_j C_j Y_j \bigg)^2 \Bigg]}\\
 & \leq \bar c \, \sqrt{\EE_{p_c} \Bigg[ \bigg( \sum_j C_j Y_j \bigg)^2 \Bigg]} = \bar c \, \sqrt{ \EE_{p_c} \Bigg[ \sum_j C_j^2 Y_j^2 \Bigg]}\\
 & \leq \bar c \, \sqrt{ \EE_{p_c} \Bigg[ \sum_j Y_j^2 \Bigg]} = \bar c \, \sqrt{ \sum_j \EE_{p_c} \big[Y_j \big]}\\
 & \leq  \hat c \, \sqrt{n^2 \pi_2(n)} = \hat c \, n \sqrt{\pi_2(n)},
\end{align*}
where the four inequalities follow, respectively, from the Cauchy-Schwarz inequality, \eqref{eq:G-EZ-bnd}, the fact that $|C_j| \leq 1$, and \eqref{eq:G-B6}, and where the first equality follows from \eqref{eq:G-Bij}, and the second one from the fact that $Y_j^2 = Y_j$.

Since the sum over $j$ of the r.h.s. of \eqref{eq:G-B8} is of order $n^2 \pi_4(3 n)$, we get that, for some universal constant $\tilde c$,
$$\pi_4(3 n) \leq \frac{\tilde c}{n} \sqrt{\pi_2(n)}.$$
This (using also that $\pi_4$ is decreasing) completes the proof of Theorem \ref{thm:main}.

\subsection*{Acknowledgments}

We thank Christophe Garban for valuable comments related to Sections~\ref{sec:Sch_St} and \ref{sec:Garban}. Many of our discussions on the subject of this paper took place at the University of Cambridge in the Fall of 2019, and we thank the Department of Pure Mathematics and Mathematical Statistics for its hospitality.

\bibliographystyle{plain}
\bibliography{Four_arm}

\end{document}